\documentclass[12pt]{amsart}
\usepackage{etex}
\usepackage{amscd,amssymb,amsmath,multicol,float,pictex}
\restylefloat{table}
\newtheorem{thm}[equation]{Theorem}
\numberwithin{equation}{section}
\newtheorem{cor}[equation]{Corollary}

\newtheorem{lem}[equation]{Lemma}

\begin{document}
\raggedbottom \voffset=-.7truein \hoffset=0truein \vsize=8truein
\hsize=6truein \textheight=8truein \textwidth=6truein
\baselineskip=18truept

\def\mapright#1{\ \smash{\mathop{\longrightarrow}\limits^{#1}}\ }
\def\mapleft#1{\smash{\mathop{\longleftarrow}\limits^{#1}}}
\def\mapup#1{\Big\uparrow\rlap{$\vcenter {\hbox {$#1$}}$}}
\def\mapdown#1{\Big\downarrow\rlap{$\vcenter {\hbox {$\ssize{#1}$}}$}}
\def\mapne#1{\nearrow\rlap{$\vcenter {\hbox {$#1$}}$}}
\def\mapse#1{\searrow\rlap{$\vcenter {\hbox {$\ssize{#1}$}}$}}
\def\mapr#1{\smash{\mathop{\rightarrow}\limits^{#1}}}
\def\ss{\smallskip}
\def\vp{v_1^{-1}\pi}
\def\at{{\widetilde\alpha}}
\def\sm{\wedge}
\def\la{\langle}
\def\ra{\rangle}
\def\on{\operatorname}
\def\spin{\on{Spin}}
\def\kbar{{\overline k}}
\def\qed{\quad\rule{8pt}{8pt}\bigskip}
\def\ssize{\scriptstyle}
\def\a{\alpha}
\def\bz{{\Bbb Z}}
\def\im{\on{im}}
\def\ct{\widetilde{C}}
\def\ext{\on{Ext}}
\def\sq{\on{Sq}}
\def\eps{\epsilon}
\def\ar#1{\stackrel {#1}{\rightarrow}}
\def\br{{\bold R}}
\def\bC{{\bold C}}
\def\bA{{\bold A}}
\def\bB{{\bold B}}
\def\bD{{\bold D}}
\def\bh{{\bold H}}
\def\bQ{{\bold Q}}
\def\bP{{\bold P}}
\def\bx{{\bold x}}
\def\bo{{\bold{bo}}}
\def\si{\sigma}
\def\Ebar{{\overline E}}
\def\dbar{{\overline d}}
\def\Sum{\sum}
\def\tfrac{\textstyle\frac}
\def\tb{\textstyle\binom}
\def\Si{\Sigma}
\def\w{\wedge}
\def\equ{\begin{equation}}
\def\b{\beta}
\def\G{\Gamma}
\def\g{\gamma}
\def\k{\kappa}
\def\psit{\widetilde{\Psi}}
\def\tht{\widetilde{\Theta}}
\def\psiu{{\underline{\Psi}}}
\def\thu{{\underline{\Theta}}}
\def\aee{A_{\text{ee}}}
\def\aeo{A_{\text{eo}}}
\def\aoo{A_{\text{oo}}}
\def\aoe{A_{\text{oe}}}
\def\fbar{{\overline f}}
\def\endeq{\end{equation}}
\def\sn{S^{2n+1}}
\def\zp{\bold Z_p}
\def\A{{\cal A}}
\def\P{{\mathcal P}}
\def\cj{{\cal J}}
\def\zt{{\bold Z}_2}
\def\bs{{\bold s}}
\def\bof{{\bold f}}
\def\bq{{\bold Q}}
\def\be{{\bold e}}
\def\Hom{\on{Hom}}
\def\ker{\on{ker}}
\def\coker{\on{coker}}
\def\da{\downarrow}
\def\colim{\operatornamewithlimits{colim}}
\def\zphat{\bz_2^\wedge}
\def\io{\iota}
\def\Om{\Omega}

\def\e{{\cal E}}
\def\exp{\on{exp}}
\def\wbar{{\overline w}}
\def\xbar{{\overline x}}
\def\ybar{{\overline y}}
\def\zbar{{\overline z}}
\def\ebar{{\overline e}}
\def\nbar{{\overline n}}
\def\rbar{{\overline r}}
\def\et{{\widetilde E}}
\def\ni{\noindent}
\def\coef{\on{coef}}
\def\den{\on{den}}
\def\lcm{\on{l.c.m.}}
\def\vi{v_1^{-1}}
\def\ot{\otimes}
\def\psibar{{\overline\psi}}
\def\mhat{{\hat m}}
\def\exc{\on{exc}}
\def\ms{\medskip}
\def\ehat{{\hat e}}
\def\etao{{\eta_{\text{od}}}}
\def\etae{{\eta_{\text{ev}}}}
\def\dirlim{\operatornamewithlimits{dirlim}}
\def\gt{\widetilde{L}}
\def\lt{\widetilde{\lambda}}
\def\st{\widetilde{s}}
\def\ft{\widetilde{f}}
\def\sgd{\on{sgd}}
\def\lfl{\lfloor}
\def\rfl{\rfloor}
\def\ord{\on{ord}}
\def\gd{{\on{gd}}}
\def\rk{{{\on{rk}}_2}}
\def\nbar{{\overline{n}}}
\def\lg{{\on{lg}}}
\def\N{{\Bbb N}}
\def\Z{{\Bbb Z}}
\def\Q{{\Bbb Q}}
\def\R{{\Bbb R}}
\def\C{{\Bbb C}}
\def\l{\left}
\def\r{\right}
\def\mo{\on{mod}}
\def\vexp{v_1^{-1}\exp}
\def\notimm{\not\subseteq}
\def\Remark{\noindent{\it  Remark}}

\def\*#1{\mathbf{#1}}
\def\0{$\*0$}
\def\1{$\*1$}
\def\22{$(\*2,\*2)$}
\def\33{$(\*3,\*3)$}
\def\ss{\smallskip}
\def\ssum{\sum\limits}
\def\dsum{{\displaystyle{\sum}}}
\def\la{\langle}
\def\ra{\rangle}
\def\on{\operatorname}
\def\o{\on{o}}
\def\U{\on{U}}
\def\lg{\on{lg}}
\def\a{\alpha}
\def\bz{{\Bbb Z}}
\def\eps{\varepsilon}
\def\br{{\bold R}}
\def\bc{{\bold C}}
\def\bN{{\bold N}}
\def\nut{\widetilde{\nu}}
\def\tfrac{\textstyle\frac}
\def\b{\beta}
\def\G{\Gamma}
\def\g{\gamma}
\def\zt{{\bold Z}_2}
\def\zth{{\bold Z}_2^\wedge}
\def\bs{{\bold s}}
\def\bx{{\bold x}}
\def\bof{{\bold f}}
\def\bq{{\bold Q}}
\def\be{{\bold e}}
\def\lline{\rule{.6in}{.6pt}}
\def\xb{{\overline x}}
\def\xbar{{\overline x}}
\def\ybar{{\overline y}}
\def\zbar{{\overline z}}
\def\ebar{{\overline \be}}
\def\nbar{{\overline n}}
\def\rbar{{\overline r}}
\def\Mbar{{\overline M}}
\def\et{{\widetilde e}}
\def\ni{\noindent}
\def\ms{\medskip}
\def\ehat{{\hat e}}
\def\xhat{{\widehat x}}
\def\nbar{{\overline{n}}}
\def\minp{\min\nolimits'}
\def\N{{\Bbb N}}
\def\Z{{\Bbb Z}}
\def\Q{{\Bbb Q}}
\def\R{{\Bbb R}}
\def\C{{\Bbb C}}
\def\S{\mathcal S}
\def\el{\ell}
\def\TC{\on{TC}}
\def\dstyle{\displaystyle}
\def\ds{\dstyle}
\def\Remark{\noindent{\it  Remark}}
\title
{Real projective space as a space of planar polygons}
\author{Donald M. Davis}
\address{Department of Mathematics, Lehigh University\\Bethlehem, PA 18015, USA}
\email{dmd1@lehigh.edu}
\date{January 15, 2015}

\keywords{Topological complexity, robotics, planar polygon spaces}
\thanks {2000 {\it Mathematics Subject Classification}: 58D29, 55R80, 70G40, 51N20
.}

\maketitle
\begin{abstract} We prove that real projective space $RP^{n-3}$ is homeomorphic to the space
$\Mbar_{n,n-2}$ of all isometry classes of $n$-gons in the plane with one side of length $n-2$ and all other sides of length 1.
This makes the topological complexity of real projective space more relevant to robotics.
 \end{abstract}

\section{Introduction}\label{intro}

The topological complexity, $\TC(X)$, of a topological space $X$ is, roughly, the number of rules required to specify how to move between
any two points of $X$.(\cite{F}) This is relevant to robotics if $X$ is the space of all configurations of a robot.

A celebrated theorem in the subject states that, for real projective space $RP^n$ with $n\ne1$, 3, or 7, $\TC(RP^n)$ is 1 greater than the dimension of the smallest
Euclidean space in which $RP^n$ can be immersed.(\cite{FY}) This is of interest to algebraic topologists because of the huge amount of work
that has been invested during the past 60 years in studying this immersion question.  See, e.g., \cite{Git}, \cite{James}, \cite{Dsur}, and \cite{Dweb}. In the popular article \cite{app}, this theorem was highlighted as an unexpected application of algebraic topology.

But, from the definition of $RP^n$, all that $\TC(RP^n)$ really tells is how hard it is to move efficiently between lines through the origin in $\R^{n+1}$, which is probably
not  very useful for robotics. Here we show that $RP^n$ may be interpreted to be the space of all polygons  of a certain type in the plane.
The edges of polygons can be thought of as linked arms of a robot, and so $\TC(RP^n)$ can be interpreted as telling how many rules are required to tell such
a robot how to move from any configuration to any other.

Let $M_{n,r}$ denote the moduli space of all oriented $n$-gons in the plane with one side of length $r$ and the rest of length 1, where two such polygons are identified if one can be obtained from the other by an orientation-preserving isometry of the plane. These $n$-gons allow  sides to intersect. Since any such $n$-gon can be uniquely rotated so that its $r$-edge is oriented in the negative $x$-direction, we can fix vertices $\bx_0=(0,0)$ and $\bx_{n-1}=(r,0)$ and define
\begin{equation}\label{Mdef}M_{n,r}=\{(\bx_1,\ldots,\bx_{n-2}): d(\bx_{i-1},\bx_i)=1,\ 1\le i\le n-1\}.\end{equation}
Here $d$ denotes distance between points in the plane.

Most of our work is devoted to proving the following theorem.
\begin{thm}\label{thm} If $n-2\le r<n-1$, then there is a $\Z/2$-equivariant homeomorphism $\Phi:M_{n,r}\to S^{n-3}$, where the involutions are reflection across the
$x$-axis in $M_{n,r}$, and the antipodal action in the sphere.\end{thm}
\noindent That this $M_{n,r}$ and $S^{n-3}$ are homeomorphic was proved earlier in \cite{HS}, but there the $\Z/2$-equivariance is not clear.

Taking the quotient of our homeomorphism by the $\Z/2$-action yields  our main new result.
It deals with the space $\Mbar_{n,r}$ of isometry classes of planar $(1^{n-1},r)$-polygons. This could be defined as the quotient of (\ref{Mdef}) modulo reflection across the $x$-axis.
\begin{cor} If $n-2\le r<n-1$, then $\Mbar_{n,r}$ is homeomorphic to $RP^{n-3}$.\end{cor}

 The fact that $\Mbar_{n,r}$ and $RP^{n-3}$ have isomorphic mod-2 cohomology rings can be extracted from earlier work of \cite{HK}.
 Noticing this is what led the author to discover the homeomorphism.

 \section{Proof of Theorem \ref{thm}}
 In this section we prove Theorem \ref{thm}.  Let $J^m$ denote the $m$-fold Cartesian product of the interval $[-1,1]$, and $S^0=\{\pm1\}$. Our model for $S^{n-3}$ is the quotient of $J^{n-3}\times S^0$
 by the relation that if any component of $J^{n-3}$ is $\pm1$, then all subsequent coordinates are irrelevant. That is, if $t_i=\pm1$, then
 \begin{equation}\label{reln}(t_1,\ldots, t_i,t_{i+1},\ldots,t_{n-2})\sim(t_1,\ldots, t_i,t'_{i+1},\ldots,t'_{n-2})\end{equation}
 for any $t'_{i+1},\ldots,t'_{n-2}$. This is just the iterated unreduced suspension of $S^0$, and the antipodal map is negation in all coordinates.
 An explicit homeomorphism of this model with the standard $S^{n-3}$ is given by
 $$(t_1,\ldots,t_{n-2})\leftrightarrow (x_1,\ldots,x_{n-2}),$$
 with $$x_i=t_i\prod_{j=1}^{i-1}\sqrt{1-t_j^2},\qquad t_i=\frac{x_i}{\sqrt{1-x_1^2-\cdots-x_{i-1}^2}}\text{ if }\sum_{j=1}^{i-1}x_i^2<1.$$
 Then $t_i=\pm1$ for the smallest $i$ for which $x_1^2+\cdots+x_i^2=1$.

Let $\P\in M_{n,r}$ be a polygon with vertices $\bx_i$ as in (\ref{Mdef}). We will define the coordinates $t_i=\phi_i(\P)$ of $\Phi(\P)$ under the homeomorphism $\Phi$ of Theorem \ref{thm}.

For $0\le i\le n-2$, we have
\begin{equation}\label{triin}n-2-i\le d(\bx_i,\bx_{n-1})\le n-1-i.\end{equation}
The first inequality follows by induction on $i$ from the triangle inequality and its validity when $i=0$. The second inequality also uses the triangle inequality together with the fact that you can get from $\bx_i$ to $\bx_{n-1}$ by $n-1-i$ unit segments. The second inequality is strict if $i=0$ and is equality if $i=n-2$.
Let $i_0$ be the minimum value of $i$ such that equality holds in this second inequality. Then the vertices $\bx_{i_0}\ldots,\bx_{n-1}$ must lie evenly spaced along a straight line segment.

Let $C(\bx,t)$ denote the circle of radius $t$ centered at $\bx$.
The inequalities (\ref{triin}) imply that, for $1\le i\le i_0$, $C(\bx_{n-1},n-1-i)$ cuts off an arc of $C(\bx_{i-1},1)$, consisting of points $\bx$ on $C(\bx_{i-1},1)$ for which $d(\bx,\bx_{n-1})\le n-1-i$. Parametrize this arc linearly, using parameter values $-1$ to 1 moving counterclockwise. The vertex $\bx_i$ lies on this arc. Set $\phi_i(\P)$ equal to the parameter value of $\bx_i$. If $i=i_0$, then $\phi_i(\P)=\pm1$, and conversely.

The following diagram illustrates a polygon with $n=7$, $r=5.2$, and $i_0=5$. We have denoted the vertices by their subscripts. The circles from left to right are $C(\bx_i,1)$ for $i$ from 0 to 4. The arcs centered at $\bx_6$ have radius 1 to 5 from right to left. We have, roughly, $\Phi(\P)=(.7,.6,.5,-.05,1)$.

  $$
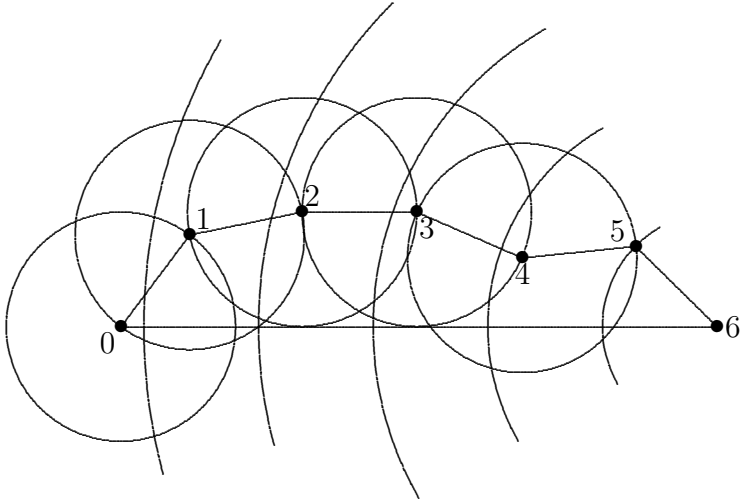
\beginpicture
\setcoordinatesystem units <.6in, .6in>
\setplotarea x from -1 to 7, y from -1.5 to 3
\circulararc 360 degrees from 1 0 center at 0 0
\circulararc 90 degrees from 4.7 .866 center at 5.2 0
\circulararc 90 degrees from 4.2 1.732 center at 5.2 0
\circulararc 90 degrees from 3.7 2.598 center at 5.2 0
\circulararc 60 degrees from 2.372 2.828 center at 5.2 0
\circulararc 45 degrees from .87 2.5 center at 5.2 0
\circulararc 360 degrees from 1.6 .8 center at .6 .8
\circulararc 360 degrees from 2.58 1 center at 1.58 1
\circulararc 360 degrees from 3.58 1 center at 2.58 1
\plot 0 0 .6 .8 1.58 1 2.58 1  3.5 .6 4.49 .7 5.2 0 0 0 /
\circulararc 360 degrees from 4.5 .6 center at 3.5 .6
\put {$\bullet$} at 5.2 0
\put {$6$} [l] at 5.27 0
\put {$\bullet$} at 0 0
\put {$0$} [tr] at -.05 -.05
\put {$\bullet$} at .6 .8
\put {${1}$} [bl] at .65 .85
\put {$\bullet$} at 1.58 1
\put {$2$} [bl] at 1.6 1.05
\put {$\bullet$} at 2.58 1
\put {$3$} [tl] at 2.6 .95
\put {$\bullet$} at 3.5 .6
\put {$4$} [t] at 3.5 .55
\put {$\bullet$} at 4.49 .7
\put {$ 5$} [br] at 4.4 .75
\endpicture$$

Here is another example, illustrating how the edges of the polygon can intersect one another, and a case with $i_0<n-2$. Again we have $n=7$ and $r=5.2$.
This time, roughly, $\Phi(\P)=(.2,-.4,.4,1,t_5)$, with $t_5$ irrelevant. Because $i_0=4$, we did not draw the circle $C(\bx_4,1)$.

  $$
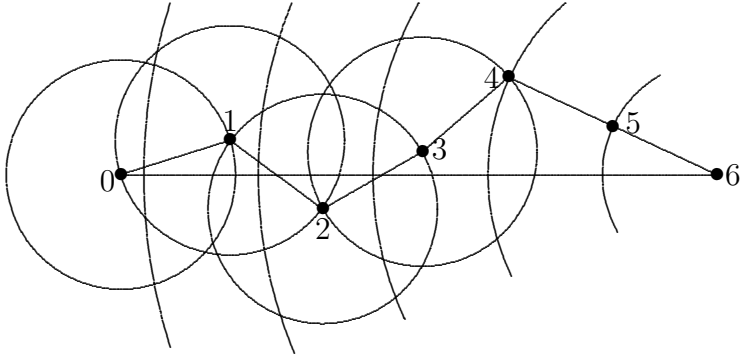
\beginpicture
\setcoordinatesystem units <.6in, .6in>
\setplotarea x from -1 to 7, y from -1.5 to 1.7
\circulararc 360 degrees from 1 0 center at 0 0
\circulararc 90 degrees from 4.7 .866 center at 5.2 0
\circulararc 75 degrees from 3.88 1.5 center at 5.2 0
\circulararc 55 degrees from 2.6 1.5 center at 5.2 0
\circulararc 45 degrees from 1.49 1.5 center at 5.2 0
\circulararc 35 degrees from .43 1.5 center at 5.2 0
\circulararc 360 degrees from 1.95 .3 center at .95 .3
\circulararc 360 degrees from 2.76 -.3 center at 1.76 -.3
\circulararc 360 degrees from 3.63 .2 center at 2.63 .2
\plot 5.2 0 0 0 .95 .3 1.76 -.3 2.63 .2 3.38 .85  5.2 0 /
\put {$\bullet$} at 5.2 0
\put {$ 6$} [l] at 5.27 0
\put {$\bullet$} at 0 0
\put {$0$} [r] at -.05 -.05
\put {$\bullet$} at .95 .3
\put {${1}$} [b] at .95 .38
\put {$\bullet$} at 1.76 -.3
\put {$2$} [t] at 1.76 -.38
\put {$\bullet$} at 2.63 .2
\put {$ 3$} [b] at 2.78 .14
\put {$\bullet$} at 3.38 .85
\put {$4$} [r] at 3.3 .85
\put {$\bullet$} at 4.29 .42
\put {$ 5$} [l] at 4.4 .455
\endpicture$$

That $\Phi$ is well defined  follows from (\ref{reln}); once we have $t_i=\pm1$, which happens first when $i=i_0$, subsequent vertices are determined and the values of subsequent
$t_j$ are irrelevant. Continuity follows from the fact that the unit circles vary continuously with the various $\bx_i$, hence so do the parameter values along the arcs cut off.
Bijectivity  follows from the construction; every set of $t_i$'s up to the first $\pm1$ corresponds to a unique polygon, and $\pm1$ will always occur. Since it maps from a compact space to a Hausdorff space, $\Phi$ is then a homeomorphism. Equivariance with respect to the involution is also clear. If you flip the polygon, you flip the whole picture, including the unit circles,
and this just negates all the $t_i$'s.

We elaborate slightly on the surjectivity of $\Phi$. The arc on $C(\bx_0,1)$ cut off by $C(\bx_{n-1},n-2)$ is determined by $n$ and $r$. Given a value of $t_1$ in $[-1,1]$, the vertex $\bx_1$ is now determined on this arc. Now the arc on $C(\bx_1,1)$ cut off by $C(x_{n-1},n-3)$ is determined, and a specified value of $t_2$ determines the vertex $\bx_2$.
All subsequent vertices of an $n$-gon are determined in this manner.

 \def\line{\rule{.6in}{.6pt}}

\end{document}